\documentclass[11pt]{article}

\usepackage{graphicx}
\usepackage{amssymb}
\usepackage{amsthm}

\usepackage{amsmath}

\usepackage{url}
\usepackage{algorithm}

\newcommand{\sgn}{\mathrm{sgn}}

\newcommand{\R}{\mathbb{R}}

\newcommand{\ve}[1]{\mbox{\boldmath $ #1$}}

\theoremstyle{definition}

\makeatletter

  \gdef\listctr{list\romannumeral\the\@listdepth}\expandafter
  
\makeatother

\newenvironment{AlgorithmSteps}[1][1]{%
  \begin{list}{\csname label\listctr\endcsname}{%
      \usecounter{\listctr}
      
      \settowidth{\labelwidth}{\textsc{Step\ #1.}}%
      \setlength{\leftmargin}{\labelwidth}\addtolength{\leftmargin}{\labelsep}}}%
  {\end{list}}

\makeatletter
\newcommand{\projalg}[1]{\ifcase #1\or GP\or SS%
    \else \@ctrerr \fi}
\makeatother

\newenvironment{ifelse}{%
  \begin{list}{}{%
      \setlength{\topsep}{0pt}\setlength{\parskip}{0pt}
      \setlength{\partopsep}{0pt}\setlength{\itemsep}{0pt}}}%
  {\end{list}}

\title{Accelerating gradient projection methods for $\ell_1$-constrained signal recovery
by steplength selection rules}

\author{I. Loris${}^1$, M. Bertero${}^2$, C. De Mol${}^3$, R. Zanella${}^4$ and L. Zanni${}^4$\\[3mm]
${}^1$Dept of Mathematics, VUB\\
${}^2$Dept of Computer and Information Sciences,
University of Genova\\
${}^3$Dept of Mathematics and ECARES, ULB\\
${}^4$Dept of Pure and Applied
Mathematics, University of Modena and Reggio Emilia}

\begin{document}

\maketitle

\begin{abstract}
We propose a new gradient projection algorithm that compares
favorably with the fastest algorithms available to date for
$\ell_1$-constrained sparse recovery from noisy data, both in
the compressed sensing and inverse problem frameworks. The
method exploits a line-search along the feasible direction and
an adaptive steplength selection based on recent strategies for
the alternation of the well-known Barzilai-Borwein rules. The
convergence of the proposed approach is discussed and a
computational study on both well-conditioned and
ill-conditioned problems is carried out for performance
evaluations in comparison with five other algorithms proposed
in the literature.
\end{abstract}

\section{Introduction}
\label{}

There has been a vast amount of recent literature dedicated to
algorithms for sparse recovery, both in the context of inverse
imaging problems and of {\em compressed sensing}. As an
alternative to the usual quadratic penalties used in
regularization theory for ill-posed or ill-conditioned inverse
problems, the use of $\ell_1$-type penalties has been advocated
in order to recover regularized solutions having sparse
expansions on a given basis or frame, such as e.g. a wavelet
system \cite{Dau04}. Denoting by $\ve x \in \R^p$ the vector of
coefficients describing the unknown object, by $\ve y \in \R^n$
the vector of (noisy) data and by $\ve K$ the linear operator
($n \times p$ matrix) modelling the link between the two, the
inverse problem amounts to finding a regularized solution of
the equation $\ve K \ve x = \ve y$. When it is known a priori
that $\ve x$ is a sparse vector, one can resort to the
following penalized least-squares strategy \cite{chendonoho},
also referred to as the {\em lasso} after Tibshirani
\cite{Tib96}:
\begin{equation}
\bar{ \ve x}(\lambda)=\arg\min_{\ve x}\| \ve K\ve x-\ve
y\|^2+2\lambda \|\ve x\|_1\label{minimizer}
\end{equation}
where $\lambda$ is a positive regularization parameter regulating
the balance between the penalty and the data misfit terms. The
norm $\| \cdot \|$ denotes the usual $\ell_2$ norm whereas $\| \ve
x\|_1 = \sum_{i=1}^p |x_i|$ is the $\ell_1$ norm of the vector
$\ve x$.

In compressed sensing (also called {\em compressive sampling}),
the aim is to reconstruct a \emph{sparse} signal or object from
a small number of linear measurements \cite{CandRomTao06,
CandTao06,Donoho2006}. The recovery of such an object can then
be achieved by searching for the sparsest solution to the
linear system $\ve K\ve x=\ve y$ representing the measurement
process, or equivalently by looking for a solution with minimum
``$\ell_0$-norm''. To avoid the combinatorial complexity of the
latter problem, one can use as a proxy a convex $\ell_1$-norm
minimization strategy. When the data $\ve y$ are affected by
measurement errors, the problem is reformulated as a penalized
least-squares optimization analogous to (\ref{minimizer}).

Let us observe that problem (\ref{minimizer}) is equivalent to the
constrained minimization problem:
\begin{equation}
\tilde {\ve x}(\rho)=\arg\min_{\| x\|_1\leq \rho}\| \ve K\ve
x-\ve y\|^2 \qquad \label{constrmin}
\end{equation}
for a certain $\rho$. One can show that $\bar {\ve x}(\lambda)$
and $\tilde {\ve x}(\rho)$ are piecewise linear functions of
$\lambda$ and $\rho$. One always has that $\bar {\ve
x}(\lambda)={\ve 0}$ for $\lambda\geq
\lambda_\mathrm{max}\equiv\max_i|({\ve K}^T{\ve y)})_i|$. The
relationship between $\lambda$ and $\rho$ is given by
$\lambda=\max_i |({\ve K}^T(\ve y-\ve K \tilde{\ve
x}(\rho)))_i|$ and $\rho=\|\bar{\ve x}(\lambda)\|_1$
\cite{DaFoL2008}.

\section{Iterative minimization algorithms}

Several iterative methods for solving the minimization problems
(\ref{minimizer}) or (\ref{constrmin}) have been proposed in
the literature. For the purpose of comparison with our new
acceleration scheme, we will focus on the following algorithms:
\begin{enumerate}
\item The Iterative Soft-Thresholding Algorithm (``ISTA'')
    proposed in
    \cite{Dau04,Figueiredo.Nowak2003,Chambolle04}) goes as
    follows: $\ve x^{(k+1)}= S_\lambda[\ve x^{(k)}+ \ve
    r^{(k)}]$ where $ \ve r^{(k)}= \ve K^T(\ve y-\ve K \ve
x^{(k)})$ is the residual in step $k$ and the (nonlinear)
soft-thresholding operator acts componentwise as $(
S_\lambda[\ve x])_i = x_i - \lambda\ \sgn(x_i)$ if
$\displaystyle{|x_i| > \lambda}$ and zero otherwise. For
any initial vector $\ve x^{(0)}$ and under the condition
$\|\ve K\|<1$, this scheme has been shown to converge to
the minimizer $\bar{\ve x}(\lambda)$ defined by
(\ref{minimizer}) \cite{Dau04}. When reinterpreted as a
forward-backward proximal scheme, convergence can be seen
to hold also for $\|\ve K\|<\sqrt{2}$
\cite{Combettes.Wajs2005}. \label{tlwalg}

\item The Fast Iterative Soft-Thresholding Algorithm
    (``FISTA''), proposed in \cite{Teboulle2008}, is a
    variation of ISTA. Defining the operator $T$ by $T(\ve
    x)= S_\lambda[\ve x+\ve K^T(\ve y-\ve K\ve x)]$, the
    FISTA algorithm is:
\begin{equation}
\ve x^{(k+1)}=T\left(\ve x^{(k)}+\frac{t^{(k)}-1}{t^{(k+1)}}
\left(\ve x^{(k)}-\ve x^{(k-1)}\right)\right), \label{fista}
\end{equation}
where ${\ve x}^{(0)}={\ve 0}$,
$t^{(k+1)}=\frac{1+\sqrt{1+4(t^{(k)})^2}}{2}$ and
$t^{(0)}=1$. It has virtually the same complexity as the
ISTA algorithm, but can be shown to have better convergence
properties.\label{fistaalg}

\item The GPSR algorithm proposed in \cite{Wright}.

\item The SPARSA algorithm proposed in \cite{sparsa}.

\item The Projected Steepest Descent (``PSD'') method
    proposed in \cite{DaFoL2008}: $\ve
    x^{(k+1)}=P_\Omega[\ve x^{(k)}+ \beta^{(k)}\ve
    r^{(k)}]$, with $\beta^{(k)}=\|\ve r^{(k)}\|^2/\|\ve K
\ve r^{(k)}\|^2$. $P_\Omega$ denotes the projection onto
the $\ell_1$-ball $\Omega$ of radius $\rho$.\label{psdalg}
\end{enumerate}
The Figures in Section \ref{sec4} provide a visual way to
compute the performance of these algorithms in two problem
examples. Note that these are the same as in \cite{Loris09},
where the reader can find comparisons to yet other methods,
including e.g. the $\ell_1$-ls method, an interior point
algorithm proposed in \cite{Kim.Koh.ea2007}.

\section{Gradient Projection with Adaptive Steplength Selection}

In this section we describe the acceleration scheme we propose
for solving the optimization problem (\ref{constrmin}). This
problem is a particular case of the general problem of
minimizing a convex and continuously differentiable function
$f(\ve x)$ over a closed convex set $\Omega\subset \R^p$. Here
$\Omega=\{\ve x \in \R^p,\ \| \ve x\|_1\leq \rho\}$. A gradient
projection method for solving this problem can be stated as in
Algorithm \ref{GPM}.

Some comments about the main steps of Algorithm GP are in
order.\\
First of all, it is worth to stress that any choice of the
steplength $\alpha_k$ in a closed interval  is permitted. This
is very important from a practical point of view since it
allows to make the updating rule of $\alpha_k$ problem-related
and oriented at optimizing the performance.\\
If the projection performed in step 2 returns a vector $\ve
h^{(k)}$ equal to $\ve x^{(k)}$, then  $\ve x^{(k)}$ is a
stationary point and the algorithm stops. When  $\ve h^{(k)}
\ne \ve x^{(k)}$, it is possible to prove that $\ve d^{(k)}$ is
a descent direction for $f$ in $\ve x^{(k)}$ and the
backtracking loop in step 5 terminates with a finite number of
runs; thus
the algorithm is well defined \cite{Bertsekas,Birgin03,BZZ}.\\
The nonmonotone line-search strategy implemented in step 5
ensures that $f(\ve x^{(k+1)})$ is lower than the maximum of
the objective function in the last $M$ iterations
\cite{GrippoLamparielloLucidi}; of course, if $M=1$ then the
strategy reduces to the standard monotone Armijo rule
\cite{Bertsekas}.

Concerning the convergence properties of the algorithm, the
following result can be derived from the analysis carried out
in \cite{Birgin03,BZZ} for more general gradient projection
schemes: if the level set $\Omega_0=\{\ve x \in \Omega : f(\ve
x) \leq f(\ve x^{(0)})\}$ is bounded, then every accumulation
point of the sequence $\{\ve x^{(k)}\}$ generated by the
Algorithm GP is a stationary point of $f(\ve x)$ in $\Omega$.
We observe that the assumption is trivially satisfied for
problem \eqref{constrmin} since in this case the feasible
region $\Omega$ is bounded.

Now, we may discuss the choice of the steplengths $\alpha_k \in
[\alpha_\mathrm{min},\alpha_\mathrm{max}]$. Steplength
selection rules in gradient methods have received an increasing
interest in the last years
 from both the theoretical and the practical point of view.
On one hand, following the original ideas of Barzilai and
Borwein (BB) \cite{BB}, several steplength updating strategies
have been devised to accelerate the slow convergence exhibited
in most cases by standard gradient methods, and a lot of effort
has been put into explaining the effects of these strategies
\cite{DaiFletcher_Asym,DaiFletcherPj,CBB,fletcher01,Frassoldati,Friedlander,DaiCG}.
On the other hand, numerical experiments on randomly generated,
library and real-life test problems have confirmed the
remarkable convergence rate improvements involved by some
BB-like steplength selections
\cite{DaiFletcherPj,CBB,Wright,Frassoldati,Thomas,Zanni,DaiCG}.
Thus, it seems natural to equip a gradient projection method
with a steplength selection that takes into account the recent
advances on the BB-like updating rules.

\begin{algorithm}[t]
\caption{(Gradient Projection Method)} \label{GPM}
%\begin{AlgorithmSteps}[4]
Choose the starting point $\ve x^{(0)}\in \Omega$, set the
parameters $\beta, \theta\in (0,1)$,
$0< \alpha_{min} <\alpha_\mathrm{max}$ and fix a positive integer $M$.\\
{\textsc{For}} $k=0,1,2,...$ \textsc{do the following steps:}
%\begin{itemize}
%\item[]
\begin{AlgorithmSteps}[4]
%\item \textsc{If} $\varphi (\ve x^k, \ve w^k) < \tau$, \textsc{then} go to Step 8, \textsc{Endif}
\item[1] Choose the parameter $\alpha_k \in
    [\alpha_\mathrm{min},\alpha_\mathrm{max}]$; \item[2]
    Projection: $\ve h^{(k)} = P_{\Omega}(\ve
    x^{(k)}-\alpha_k\nabla f(\ve
x^{(k)}))$;\\
If $\ve h^{(k)}=\ve x^{(k)}$ then stop, declaring that $\ve
x^{(k)}$ is a stationary point;
 \item[3] Descent direction: $\ve d^{(k)} = \ve h^{(k)}-\ve
     x^{(k)}$; \item[4] Set $\lambda_k = 1$ and
     $f_\mathrm{max}= \max_{0\leq j \leq \min(k,M-1)}f(\ve
     x^{(k-j)})$;
\item[5] Backtracking loop:
\begin{ifelse} \item \textsc{If} $f(\ve x^{(k)}+\lambda_k \ve
d^{(k)})\leq f_\mathrm{max}+\beta\lambda_k\nabla f(\ve
           x^{(k)})^T \ve d^{(k)}$ \ \textsc{then} \\
       \hspace*{.8cm}  go to Step 6;
\item \textsc{Else} \\ \hspace*{.8cm} set $\lambda_k = \theta
    \lambda_k$ and go to Step 5; \item \textsc{Endif}
\end{ifelse}
\item[6] Set $\ve x^{(k+1)} = \ve x^{(k)} + \lambda_k  \ve
    d^{(k)}$.
\end{AlgorithmSteps}
%\end{itemize}
\textsc{End}
%\end{AlgorithmSteps}
\end{algorithm}

\noindent First of all we must recall the two BB rules usually
exploited by the main steplength updating strategies. To this
end, by denoting with $I$ the $p\times p$ identity matrix, we
can regard the matrix $B(\alpha_k)=(\alpha_k I)^{-1}$ as an
approximation of the Hessian $\nabla^2 f(\ve{x}^{(k)})$ and
derive two updating rules for $\alpha_k$ by forcing
quasi-Newton properties on $B(\alpha_k)$:
\begin{equation}
\alpha_k^{\text{BB1}} = \mbox{arg}\!\min_{\alpha_k\in\mathbb{R}}
\|B(\alpha_k)\ve{s}^{(k-1)} - \ve{z}^{(k-1)} \|
\quad \mathrm{and}\quad
\alpha_k^{\text{BB2}} = \mbox{arg}\!\min_{\alpha_k\in\mathbb{R}}
\|\ve{s}^{(k-1)} - B(\alpha_k)^{-1}\ve{z}^{(k-1)} \|,
\end{equation}
where $\ve{s}^{(k-1)} = \ve{x}^{(k)} - \ve{x}^{(k-1)}$ and
$\ve{z}^{(k-1)} = \nabla f(\ve{x}^{(k)}) - \nabla
f(\ve{x}^{(k-1)})$. In this way, the steplengths
\begin{equation}
\alpha_k^{\text{BB1}} = \frac{{\ve{s}^{(k-1)}}^T
\ve{s}^{(k-1)}}{{\ve{s}^{(k-1)}}^T  \ve{z}^{(k-1)}}\,, \qquad
\qquad \alpha_k^{\text{BB2}} = \frac{{\ve{s}^{(k-1)}}^T
\ve{z}^{(k-1)} } {{\ve{z}^{(k-1)}}^T \ve{z}^{(k-1)}}\,,
\label{bbrules}
\end{equation}
are obtained.

\begin{algorithm}[t]
\caption{(Steplength Selection for GP)} \label{SS} \textsc{if}
\ $k=0$ \ \textsc{then}
\\ \hspace*{.5cm}
%arbitrarily
set $\alpha_0 \in [\alpha_\mathrm{min},\alpha_\mathrm{max}]$, $\tau_1\in (0,1)$ and a non-negative integer $M_{\alpha}$;\\[.1cm]
\textsc{else}\\[.1cm]
%Let $k>0$ and consider the iteration $k$ of the SGP. \\[.2cm]
\hspace*{.5cm}\textsc{if} ${\ve{s}^{(k-1)}}^T \ve{z}^{(k-1)} \le 0$ \  \textsc{then}\\
\hspace*{1cm}$\alpha_k = \alpha_\mathrm{max}$;\\
\hspace*{.5cm}\textsc{else}\\
\hspace*{1cm}$\alpha_k^{(1)} = \max\left\{\alpha_\mathrm{min},
\min\left\{\frac{{\ve{s}^{(k-1)}}^T
\ve{s}^{(k-1)}}{{\ve{s}^{(k-1)}}^T \ve{z}^{(k-1)} },
 \alpha_\mathrm{max}\right\}\right\}$;\\[1mm]
\hspace*{1cm}$\alpha_k^{(2)} = \max\left\{\alpha_\mathrm{min},
\min\left\{\frac{{\ve{s}^{(k-1)}}^T  \ve{z}^{(k-1)}
}{{\ve{z}^{(k-1)}}^T  \ve{z}^{(k-1)}},
 \alpha_\mathrm{max}\right\}\right\}$;\\
%%\hspace*{.5cm}\textsc{endif}
\\[.2cm]
\hspace*{.6cm}\hspace*{.5cm}\textsc{if} $ {\alpha_k^{(2)}}/{\alpha_k^{(1)}} \le \tau_k$\ \textsc{then}\\
\hspace*{.6cm}\hspace*{1cm}$\alpha_k =
\min\left\{\alpha_j^{(2)}, \
j=\max\left\{1,k-M_{\alpha}\right\},\dots,k\right\}$;\\
\hspace*{.6cm}\hspace*{1cm}$\tau_{k+1} = \tau_{k}*0.9$;\\
\hspace*{.6cm}\hspace*{.5cm}\textsc{else}\\
\hspace*{.6cm}\hspace*{1cm}$\alpha_k = \alpha_k^{(1)}$;\\
\hspace*{.6cm}\hspace*{1cm}$\tau_{k+1} = \tau_{k}*1.1$;\\
\hspace*{.6cm}\hspace*{.5cm}\textsc{endif}\\[.2cm]
\hspace*{.6cm}\textsc{endif}\\[.1cm]
\textsc{endif}
\end{algorithm}

\noindent At this point, inspired by the steplength
alternations successfully implemented in recent gradient
methods \cite{Frassoldati,DaiCG}, we propose a steplength
updating rule for GP which adaptively alternates the values
provided by \eqref{bbrules}. The details of the GP steplength
selection are given in Algorithm \ref{SS}. This rule decides
the alternation between two different selection strategies
%$\alpha_k^{(1)}$ and $\alpha_k^{(2)}$
by means of the variable threshold $\tau_k$ instead of a
constant parameter as done in \cite{Frassoldati} and
\cite{DaiCG}. This trick makes the choice of $\tau_0$ less
important for the GP performance and, in our experience, seems
able to avoid the drawbacks due to the use of the same
steplength rule in too many consecutive iterations. In the
following we denote by GPSS the algorithm \ref{GPM} equipped
with the steplength selection \ref{SS}.
\\
We end this section by describing the setting for the GPSS parameters used in the computational study of this work:
\begin{itemize}
\item {\it line-search parameters}: $M=1$ (monotone line-search), \ \ $\theta=0.5$, \ \ $\beta=10^{-4}$;
\item {\it steplength parameters}:
    $\alpha_\mathrm{min}=10^{-10}$, \ \
    $\alpha_\mathrm{max}=10^{10}$,\\
    $\alpha_0 =
    \max\{\alpha_\mathrm{min},\min\{\| P_{\Omega}(\ve
    x^{(0)}-\nabla f(\ve
    x^{(0)}))\|_{\infty}^{-1},\alpha_\mathrm{max}\}\}$, \ \
$\tau_1=0.5$,\ \
$M_{\alpha}=2$.
\end{itemize}
In our experience the above setting often provides satisfactory
performance; however, it can not be considered optimal for
every application and a careful parameter tuning is always
advisable.

\section{Numerical experiments}

\label{sec4}

To assess the performances of our GPSS algorithm and estimate
the gain in speed it can provide with respect to the algorithms
1 to 5, we perform some numerical tests. To this purpose we
adopt the methodology proposed in \cite{Loris09} and based on
the notion of \emph{approximation isochrones}. It improves on
the comparisons made for a single value of $\lambda$ or $\rho$,
i.e. for a single level of sparsity of the recovered object.

For values of $\lambda$ in a given interval
$\lambda_\mathrm{min} \le \lambda \le \lambda_\mathrm{max}$,
one computes the minimizer $\bar {\ve x}(\lambda)$ of
(\ref{minimizer}). When the number of nonzero components in
$\bar {\ve x}(\lambda)$ is not too large, this can be done by
means of the direct (non-iterative) {\em homotopy} method
\cite{Osborne.Presnell.ea2000} or LARS algorithm
\cite{Efron.Hastie.ea2004}. Then, for a fixed and given
computation time, one runs one of the algorithms for each value
of $\lambda$ (or $\rho$). The relative error $e^{(k)}= \|\ve
x^{(k)}(\lambda)-\bar {\ve x}(\lambda)\| / \|\bar {\ve
x}(\lambda)\|$ reached at the end of the computation is plotted
as a function of $\lambda$ and hence this plot is just the
approximation isochrone showing the degree of accuracy reached
in the given amount of computing time for each value of
$\lambda$. A set of such plots allow to quickly grasp the
performances of a given algorithm in various parameter regimes
and to easily compare it with other methods; it reveals in one
glance under which circumstances the algorithms do well or
fail. The paper \cite{Loris09} also demonstrates the fact that
the relative performances of the algorithms may strongly depend
on the specific application one considers, and in particular on
the properties of the linear operator $\ve K$ modelling the
problem.

We test the different algorithms on two different operators
arising typically either from a compressed sensing or from an
inverse problem. In both cases the matrix $\ve K$ is of size
1848x8192. In the first case, the elements of $\ve K$ are taken
from a Gaussian distribution with zero mean and variance such
that $\|{\ve K}\|=1$. This matrix is rather well conditioned
and can serve as a paradigm of compressed sensing applications.
It is applied to a sparse vector and perturbed by additive
gaussian noise (about $2\%$) to yield the data $\ve y$. The
second matrix models a severely ill-conditioned linear inverse
problem that finds its origin in a problem of seismic
tomography described in detail in \cite{Loris.Nolet.ea2007}.

%gaussian matrix
\begin{figure}
\centering
\resizebox{\textwidth}{!}{\includegraphics{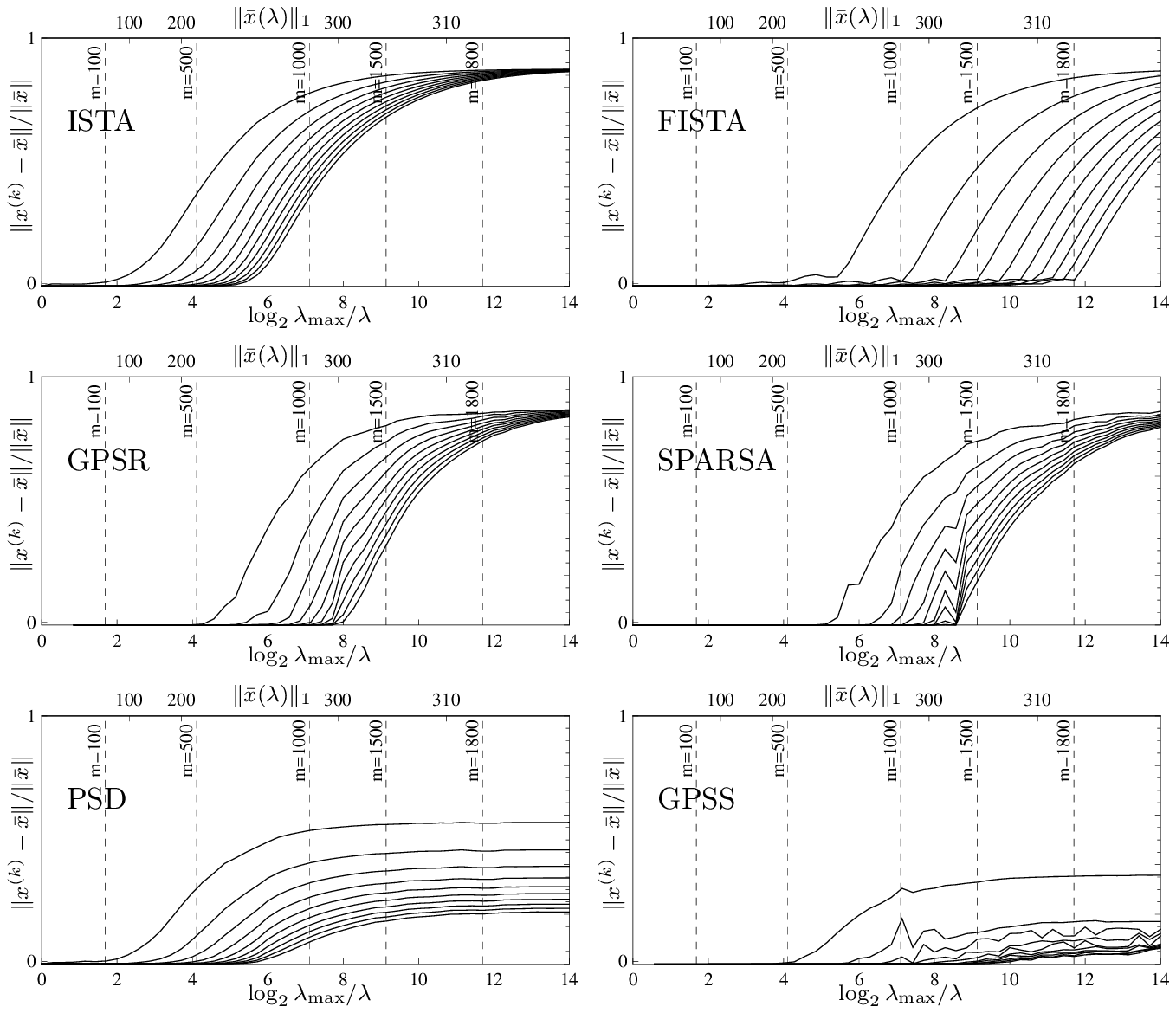}}
\caption{Approximation isochrones in the case of the gaussian
random matrix  for $t=6,12,\ldots, 60$ seconds.}\label{gausspic}
\end{figure}

For both operators, the minimizer $\bar {\ve x}(\lambda)$ is
computed for 50 different values of $\lambda$ (or equivalently,
50 different values of $\rho$). Then, for each iterative
algorithm, we make plots having the relative error $e^{(k)}$ on
the vertical axis and $\log_2\lambda_\mathrm{max}/\lambda$ on
the bottom horizontal axis (on the top horizontal axis the
value of $\rho=\|\bar{\ve x}\|_1$ is also reported). The number
of nonzero components $m$ in $\bar {\ve x}(\lambda)$ is
indicated by vertical dashed lines. In each plot we report the
isochrone lines that correspond to a given amount of computer
time. In this way one can see how close, for the different
values of $\lambda$, the iterates approach the minimizer after
a given time. Let us remark that although the reported
computing times are of course specific to a given computer and
implementation, the overall behavior of the isochrones should
be fairly general. For example, the fact that they get very
close to each other in some places can be interpreted as a
bottleneck feature of the algorithm.

\begin{figure}
\centering
\resizebox{\textwidth}{!}{\includegraphics{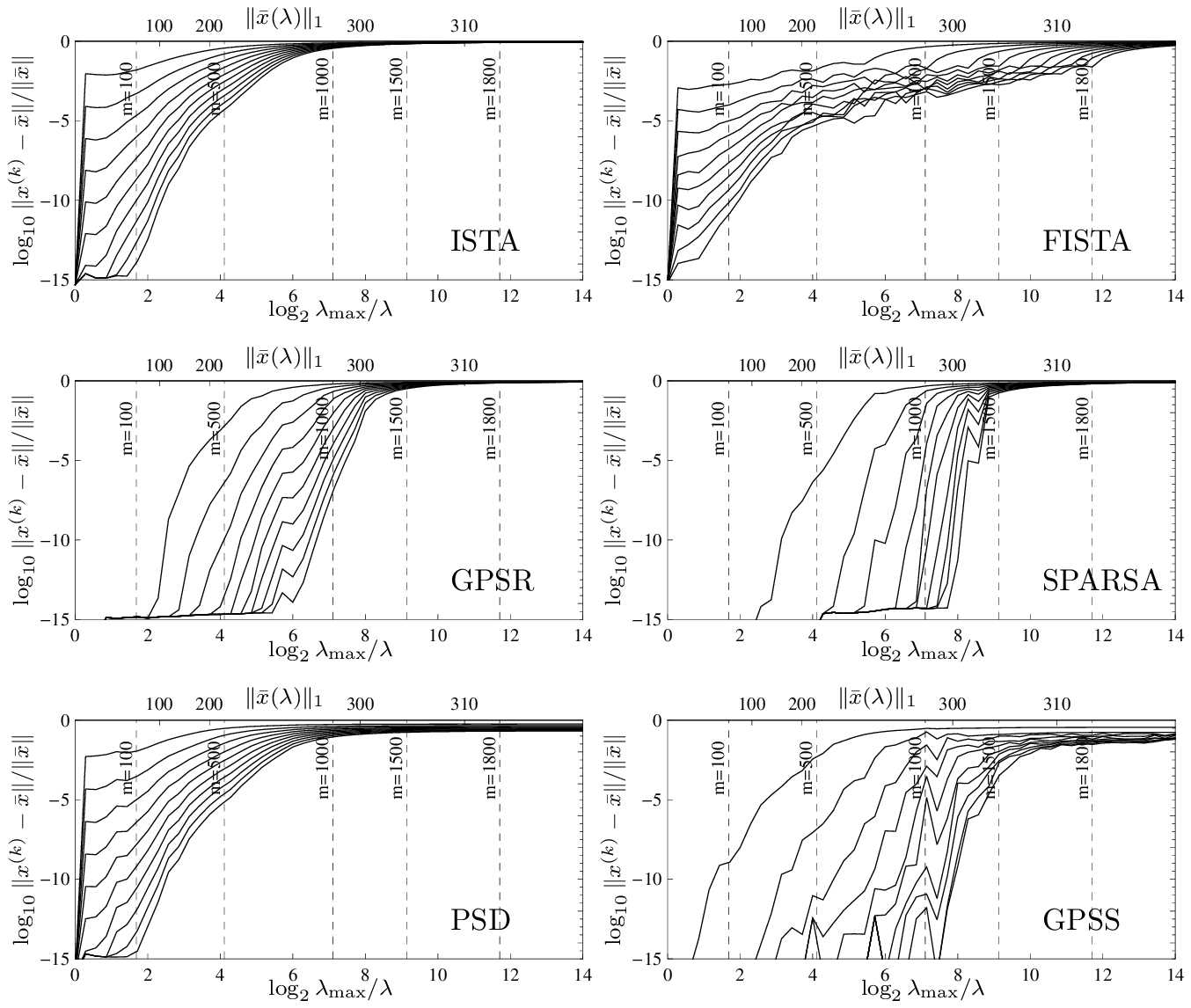}}
\caption{The same as Figure \ref{gausspic} but in a semi-log plot.}\label{gausspiclog}
\end{figure}

In Figure \ref{gausspic}, we report the results for the ISTA,
FISTA, GPSR, SPARSA, PSD and our new algorithm GPSS for the
case of the gaussian random matrix. The proposed GPSS algorithm
compares favorably with the other five, especially for small
values of $\lambda$. Experiments made by varying the parameter
$M$ showing no significant difference, we report here only the
results obtained with $M=1$ (monotonic line search). However,
the behavior for large penalties is not clearly visible on
Figure \ref{gausspic}.  It is better demonstrated when using a
logarithmic scale for the relative error on the vertical axes
as reported in Figure \ref{gausspiclog}.

%gaussian matrix log%\textwidth 9.3

% geo
\begin{figure}[h]
\centering
\resizebox{\textwidth}{!}{\includegraphics{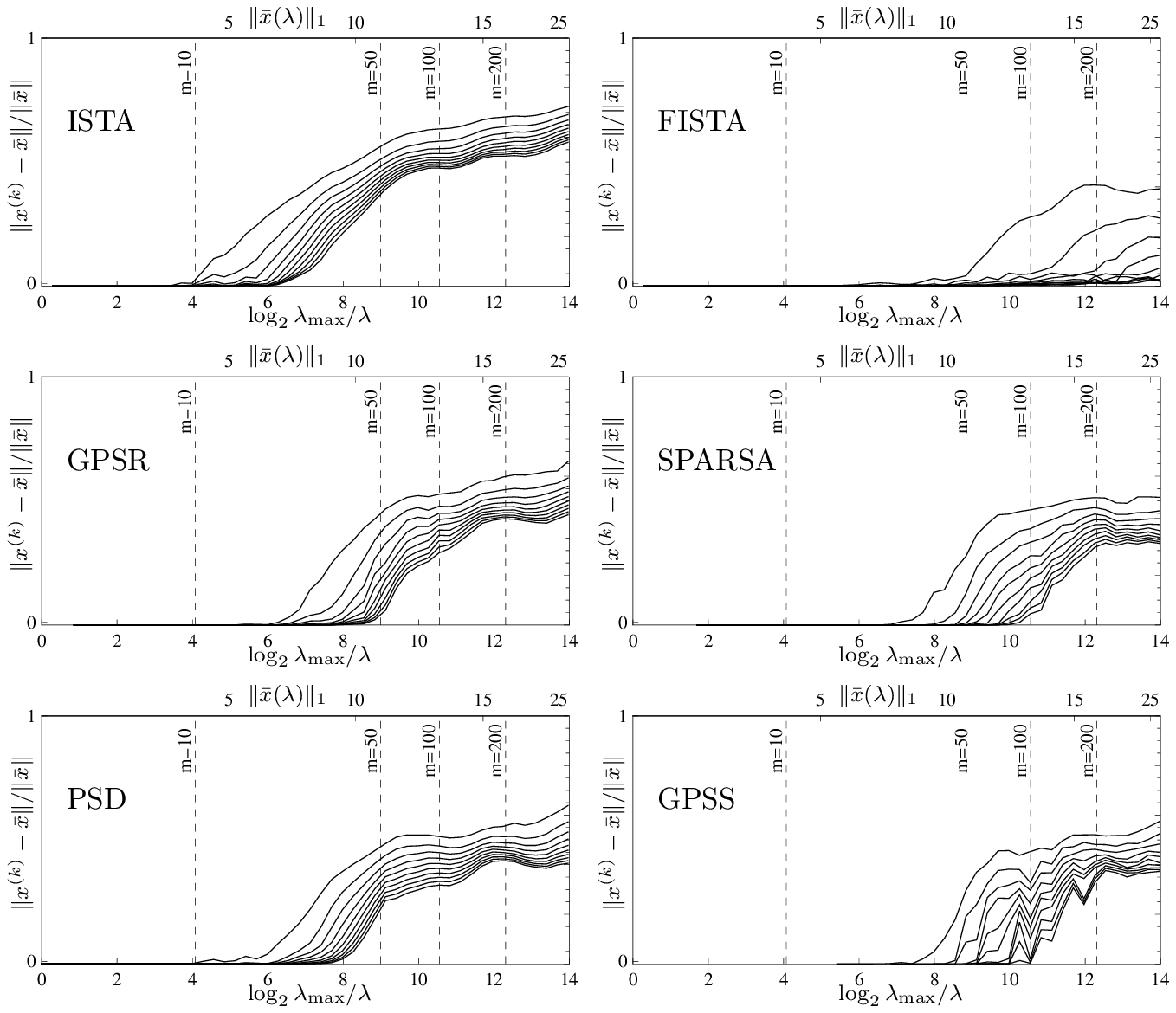}}
\caption{Approximation isochrones for the seismic inverse problem for $t=1,\ldots,10$ minutes.}\label{geopic}
\end{figure}

\begin{figure}[h]
\centering
\resizebox{\textwidth}{!}{\includegraphics{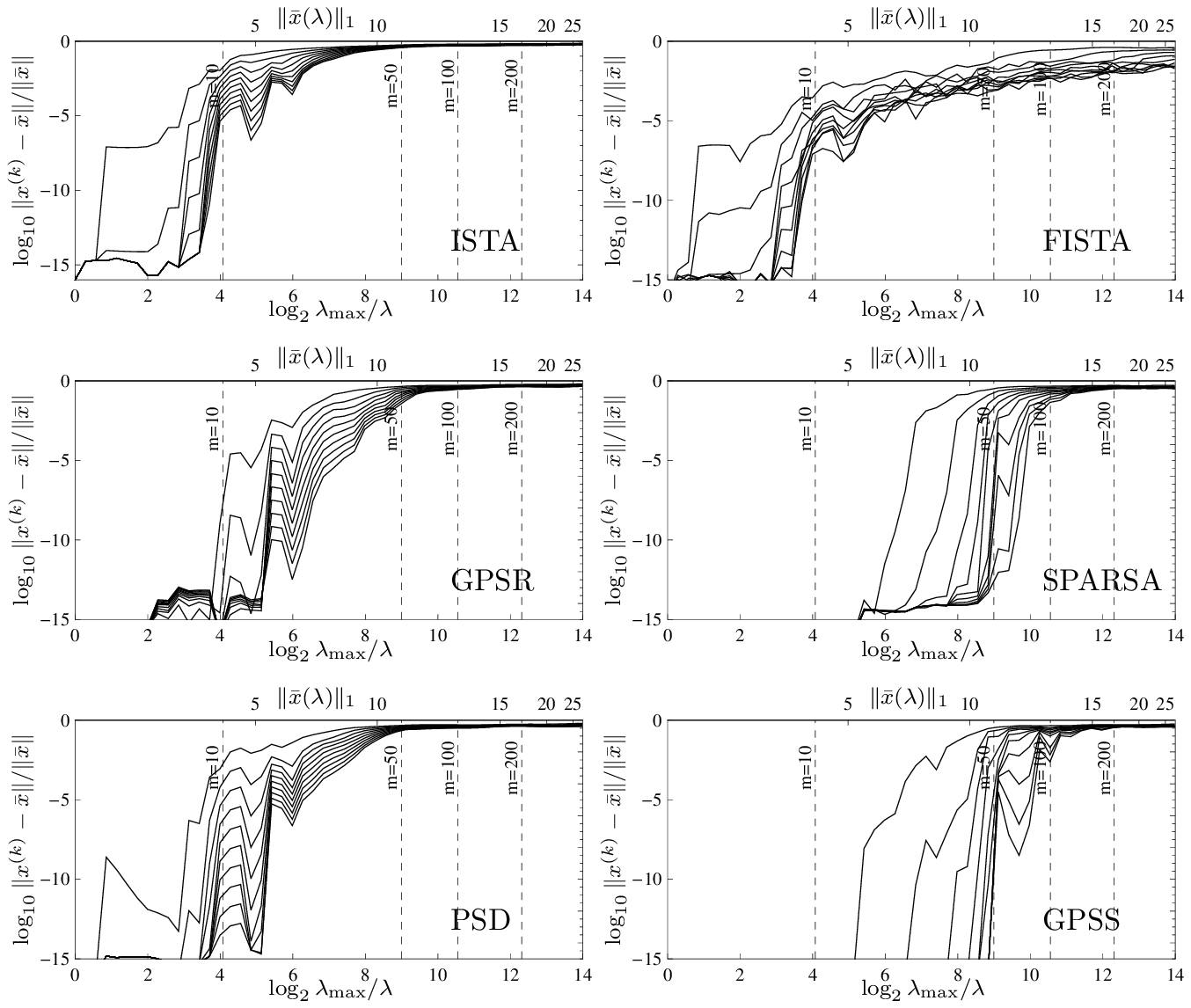}}
\caption{The same as in Figure \ref{geopic} in a semi-log plot.}\label{geopiclog}
\end{figure}

In Figures \ref{geopic} and \ref{geopiclog}, we report the
results for the case of the ill-conditioned matrix arising from
the seismic inverse problem. Clearly, for this operator, ISTA,
GPSR and PSD have a lot of difficulty in approaching the
minimizer for small values of $\lambda$ (lines not approaching
$e=0$). The FISTA algorithm appears to work best for small
penalty parameters whereas GPSS and SPARSA compete for the
second place in such instance. From Figure \ref{geopiclog}, we
see that the GPSS and SPARSA algorithms are performing best for
large values of $\lambda$.

The reported encouraging numerical results call of course for
further experiments, but we believe that they are sufficiently
representative to allow honest extrapolation to reliable
conclusions holding more generally. As seen, the proposed GPSS
algorithm performs well for the compressed sensing problem: for
small values of $\lambda$, it clearly outperforms the other
algorithms (see Figure \ref{gausspic}) whereas it is still
competitive for larger values of $\lambda$. In the
ill-conditioned inversion problem, GPSS an SPARSA appear to
perform better than all other tested algorithms for large
values of $\lambda$, whereas they are challenged by the FISTA
method for smaller values.

\section*{Acknowledgements}

I.L. and C.D.M. are supported by grant GOA-062 of the VUB. I.L.
is supported by grant G.0564.09N of the FWO-Vlaanderen. M.B.,
R.Z. and L.Z. are partly supported by MUR grant 2006018748.


\begin{thebibliography}{99}

\bibitem{BB}
J. Barzilai, J.M. Borwein, Two point step size gradient methods,
IMA J. Numer. Anal. 8 (1988) 141--148.

\bibitem{Teboulle2008}
A. Beck, M. Teboulle, A Fast Iterative Shrinkage-Thresholding
Algorithm for Linear Inverse Problems, SIAM J. Imaging Sciences,
forthcoming.

\bibitem{Bertsekas}
D.P. Bertsekas, Nonlinear Programming, Athena Scientific, 2nd
edition, 1999.


\bibitem{Birgin03} E.G. Birgin, J. M. Mart\'inez, M. Raydan,
    {Inexact spectral projected gradient methods on convex
    sets}, IMA J. Numer. Anal. {23} (2003) 539--559.

\bibitem{BZZ} S. Bonettini, R. Zanella, L. Zanni, A scaled
    gradient projection method for constrained image
    deblurring, Inverse Problems 25:015002 (2009).

\bibitem{CandRomTao06}
E. Cand\`es, J. Romberg, T. Tao, Robust uncertainty principles:
Exact signal reconstruction from highly incomplete frequency
information, IEEE Trans. Inform. Theory 52 (2006) 489--509.

\bibitem{CandTao06} E. Cand\`es,  T. Tao, Near optimal signal
    recovery from random projections: Universal encoding
    strategies?, IEEE Trans. Inform. Theory 52 (2006)
    5406--5425.

\bibitem{Chambolle04}
A. Chambolle, An Algorithm for Total Variation Minimization and
Applications, J. Math. Imaging and Vision, 20 (2004) 89--97.

\bibitem{chendonoho} S. S. Chen, D. Donoho, M. A.
    Saunders, Atomic Decomposition by Basis Pursuit, SIAM J. Sci.
    Comput., 20 (1998) 33--61.


\bibitem{Combettes.Wajs2005}
P.L. Combettes, V.R. Wajs, Signal Recovery by Proximal
Forward-Backward Splitting, Multiscale Model. Simul. 4 (2005)
1168--1200.

\bibitem{DaiFletcher_Asym}
Y.H. Dai, R. Fletcher, On the asymptotic behaviour of some new
gradient methods, Math. Programming 103 (2005) 541--559.

\bibitem{DaiFletcherPj}
Y.H. Dai, R. Fletcher, New algorithms for singly linearly
constrained quadratic programming problems subject to lower and
upper bounds, Math. Programming 106 (2006) 403--421.

\bibitem{CBB}
Y.H. Dai, W.W. Hager, K. Schittkowski, H. Zhang, The cyclic
Barzilai-Borwein method for unconstrained optimization, IMA J.
Numer. Anal. 26 (2006) 604--627.

\bibitem{Dau04}
I. Daubechies, M.~Defrise, C.~De~Mol, An iterative thresholding
algorithm for linear inverse problems with a sparsity constraint,
Comm. Pure Appl. Math. 57 (2004) 1413--1457.

\bibitem{DaFoL2008}
I. Daubechies, M. Fornasier, I. Loris, Accelerated projected
gradient method for linear inverse problems with sparsity
constraints, J. Fourier Anal. Appl. 14 (2008) 764--792.

\bibitem{Donoho2006}
D.L. Donoho, Compressed sensing, IEEE Trans. Inform. Theory 52
(2006) 1289--1306.

\bibitem{Efron.Hastie.ea2004}
B. Efron, T. Hastie, I. Johnstone, R. Tibshirani, Least angle
regression, Ann. Statist. 32 (2004) 407--499.

\bibitem{fletcher01} R. Fletcher, On the Barzilai-Borwein method,
Technical Report NA/207, Department of Mathematics, University of
Dundee, Dundee, UK, 2001.

\bibitem{Figueiredo.Nowak2003} M.A.T. Figueiredo, R.D. Nowak,
    An EM algorithm for wavelet-based image restoration, IEEE
    Trans. Image Process. 12 (2003) 906--916.

\bibitem{Wright}
M.A.T. Figueiredo, R.D. Nowak, S.J. Wright, Gradient projection
for sparse reconstruction: Application to compressed sensing and
other inverse problems, IEEE J. Selected Topics in Signal Process.
1 (2007) 586--597.

\bibitem{Frassoldati}
G. Frassoldati, G. Zanghirati, L. Zanni, New adaptive stepsize
selections in gradient methods, J. Industrial and Management
Optim. 4 (2008) 299--312.

\bibitem{Friedlander}
A. Friedlander, J.M. Mart\'inez, B. Molina, M. Raydan, Gradient
method with retards and generalizations, SIAM J. Numer. Anal. 36
(1999) 275--289.

\bibitem{GrippoLamparielloLucidi}
L. Grippo, F. Lampariello, S. Lucidi, A nonmonotone line-search
technique for Newton's method , SIAM J. Numer. Anal. 23 (1986)
707--716.

\bibitem{Kim.Koh.ea2007} S.-J. Kim, K. Koh, M. Lustig, S. Boyd,
    D. Gorinevsky, A method for
large-scale $\ell_1$-regularized least squares, IEEE Trans. on
Selected Topics in Signal Processing 1 (2007)  606--617.

\bibitem{Loris09} I. Loris, On the performance of algorithms
    for the minimization of $\ell_1$-penalized functionals,
    Inverse Problems 25:035008 (2009).

\bibitem{Loris.Nolet.ea2007}
I. Loris, G. Nolet, I. Daubechies, F.A. Dahlen, Tomographic
inversion using $\ell_1$-norm regularization of wavelet
coefficients, Geophysical Journal International 170 (2007)
359--370.

\bibitem{Osborne.Presnell.ea2000}
M.R. Osborne, B. Presnell, B.A. Turlach, A new approach to
variable selection in least squares problems, IMA J. Numer. Anal.
20 (2000) 389--403.

\bibitem{Thomas}
T. Serafini, G. Zanghirati, L. Zanni, Gradient projection methods
for quadratic programs and applications in training support vector
machines, Optim. Meth. Soft. 20 (2005) 343--378.

\bibitem{Tib96} R. Tibshirani, Regression selection and
    shrinkage via the lasso, J. R. Stat. Soc. Ser. B 58 (1996)
    267--288.

\bibitem{sparsa} S. Wright, R. Nowak, M. Figueiredo,
    Sparse reconstruction by separable approximation, IEEE Transactions on Signal
    Processing (2009), forthcoming.

\bibitem{Zanni}
L. Zanni, An improved gradient projection-based decomposition
technique for support vector machines, Comput. Management Sci. 3
(2006) 131--145.

\bibitem{DaiCG} B. Zhou, L. Gao, Y.H. Dai, Gradient methods
    with adaptive step-sizes, Comput. Optim. Appl. 35 (2006)
    69--86.

\end{thebibliography}
\end{document}